 \documentclass[centertags]{amsart}

   \newtheorem{lemma}{Lemma}[section]
   \newtheorem{theorem}[lemma]{Theorem}
   \newtheorem{remark}[lemma]{Remark}

   \newtheorem{definition}[lemma]{Definition}
   \newcommand{\eps}{\varepsilon}

\title[]{}

\author{Jinqiao Duan}
\address[J. Duan]
{Department of Applied Mathematics\\
 Illinois Institute of Technology\\
   Chicago, IL 60616, USA and Department of Mathematics\\
University of Science and Technology of China\\
Hefei, Anhui 230026,
China }
\email[J.~Duan]{duan@iit.edu}

\author{Kening Lu}
\address[K. Lu]
{Department of Mathematics\\
    Brigham Young University\\
    Provo, Utah 84602, USA
    and
    Department of Mathematics \\
    Michigan State University\\
    East Lansing, MI 48824, USA}
\email[K.~Lu]{klu@math.byu.edu or klu@math.msu.edu}

\author{Bj{\"o}rn Schmalfu{\ss}}
\address[B. Schmalfu{\ss}]{%
  Department of Sciences\\
  University of Applied Sciences\\
  Geusaer Stra{\ss}e\\
  06217 Merseburg, Germany}
\email[B.~Schmalfu{\ss}]{bjoern.schmalfuss@in.fh-merseburg.de}

\title[Invariant manifolds for SPDE's] {Invariant manifolds for stochastic partial
 differential equations}

 \date{October 21, 2002 }
\subjclass{2000 MSC -- Primary: 60H15; Secondary: 37H10, 37L55,
37L25, 37D10}
\keywords{Invariant manifolds, cocycles, non-autonomous dynamical
systems, stochastic partial differential equations, generalized
fixed points\\
This work was partially supported by NSF 0209326, NSF0200961, and
a travel grant from German DFG Schwerpunktprogramm {\em
Interagierende zuf\"allige Systeme von hoher Komplexit\"at}.}

\begin{document}

\begin{abstract}
{\bf Annals of Probability 31(2003), 2109-2135.}
 Invariant manifolds provide the geometric structures for
describing and understanding dynamics of nonlinear  systems. The
theory of invariant manifolds
 for both finite and infinite dimensional  autonomous  deterministic
 systems,
and for  stochastic ordinary differential equations
 is relatively mature. In this paper,
we present a unified theory of invariant  manifolds for infinite
dimensional {\em random} dynamical systems generated by
    {\em stochastic} partial differential equations.
We  first introduce a random graph transform and a fixed point
theorem for non-autonomous systems. Then we show the existence of
generalized fixed points which give the desired invariant
manifolds.

 \bigskip
 {\bf Short Title:}  Invariant Manifolds for SPDEs

\end{abstract}

\maketitle

%%%%   Section 1   %%%%%%%%%%
%%%%%%%%%%%%%%%%%%%%%
\section{Introduction}

Invariant manifolds are essential for describing and understanding
dynamical behavior of nonlinear and random  systems. Stable,
unstable and center manifolds have been widely used in the
investigation of  infinite dimensional {\em deterministic}
dynamical systems. In this paper, we are concerned with invariant
manifolds for {\em stochastic}  partial differential equations.

\vskip0.1in
 The theory of invariant manifolds for deterministic
dynamical systems has a long and rich history. It was first
studied by Hadamard \cite{Had01}, then, by Liapunov \cite{Lia47}
and Perron \cite{Per28} using a different approach. Hadamard's
graph transform method is a geometric approach, while
Liapunov-Perron  method is analytic in nature. Since then, there
is an extensive literature on the stable, unstable, center,
center-stable, and center-unstable manifolds for both finite and
infinite dimensional deterministic {\em autonomous} dynamical
systems (see Babin and Vishik \cite{BabVis89} or Bates et al.
\cite{BatLuZen98} and the references therein). The theory of
invariant manifolds for {\em non-autonomous} abstract semilinear
parabolic equations may be found in Henry \cite{Hen81}. Invariant
manifolds with invariant foliations for more general infinite
dimensional {\em non-autonomous} dynamical systems was studied  in
Chow et al.\cite{ChoLuLin91}. Center manifolds for infinite
dimensional {\em non-autonomous} dynamical systems was considered
in Chicone and Latushkin \cite{ChiLat97}.

 \vskip0.1in
Recently, there are some works on invariant manifolds for
stochastic or random ordinary differential equations (finite dimensional systems)
 by Wanner
\cite{Wanner}, Arnold \cite{Arn98}, Mohammed and Scheutzow
\cite{MohScheu99}, and Schmalfu{\ss} \cite{Schm97a}. Wanner's
method is based  on the Banach fixed point theorem on some Banach
space containing functions with particular exponential growth
conditions, which is essentially  the Liapunov-Perron approach. A
similar technique has been used by Arnold. In contrast to this
method, Mohammed and Scheutzow have applied a classical technique
due to Ruelle \cite{Rue82}  to  stochastic differential equations
driven by semimartingals.   In Caraballo et al. \cite{CarLanRob01}
an invariant manifold for a stochastic reaction diffusion equation
of pitchfork type has been considered. This manifold connects
different stationary solutions of the stochastic differential
equation. In Koksch and Siegmund \cite{KokSie01} the pullback
convergence has been used to construct an inertial
manifold for non-autonomous dynamical systems. \\

In this paper, we will prove the existence of an invariant
manifold for a nonlinear stochastic evolution equation with a
multiplicative white noise:
\begin{equation} \label{spde0}
\frac{d\phi}{dt}=A\phi+F(\phi)+\phi\,\dot{W}
\end{equation}
where $A$ is a generator of a $C_0$-semigroup satisfying a
exponential dichotomy condition, $F(\phi)$ is a Lipschitz
continuous operator with $F(0)=0$, and $\phi\,\dot{W}$ is the
noise. The precise conditions on them will be given in the next
section.  Some physical systems or fluid systems with noisy
perturbations proportional to the state of the system may be modeled
by this equation. \\

A similar object, inertial manifolds, has been considered by
Bensoussan and Flandoli \cite{BenFla95}, Chueshov and Girya
\cite{GirChu95} or Da Prato and Debussche \cite{DaPDeb96} for the
equations with pure white noises. Their approaches \cite{DaPDeb96}
and \cite{GirChu95}  are based on properties of Ito stochastic
differential equations like Ito's formula, martingales and Ito
integrals.\\

Here, we consider the stochastic partial differential equations
with multiplicative noises and our method is based on the theory
of random dynamical systems. In particular, we are able to
formulate conditions such that a general random evolution equation
\begin{equation}\label{spde00}
\frac{d\phi}{dt}=A\phi+G(\theta_t\omega,\phi)
\end{equation}
has an invariant manifold providing a condition on the spectral
gap and the Lipschitz constant of $G$ in $\phi$. The random
dynamical systems generated by (\ref{spde0}) and (\ref{spde00})
are conjugated, which allows us to determinate the manifold for
(\ref{spde0}) by the manifold for (\ref{spde00}).\\

Our method showing the existence of an invariant manifold is
different from the methods mentioned above, which is an extension
of the result by Schmalfu{\ss} \cite{Schm97c}.\\
However, this article only deals with a finite dimensional
equation which is semi-coupled. We will introduce a {\em random}
graph transform. In contrast to \cite{BenFla95}, \cite{GirChu95},
and \cite{DaPDeb96} this  graph transform defines a new and
lifted random dynamical system on the space of appropriate graphs.
One ingredient of a random dynamical system is a cocycle (see the
next section). An invariant graph of this graph transform is a
generalized fixed point for a cocycle. This generalized fixed
point defines an entire trajectory for the cocycle. Applying this
fixed point theorem to the graph transform dynamical system, we
can find under a gap condition a fixed point contained in the set
of Lipschitz continuous graphs which represent the invariant
manifold.\\

The main assumption is the gap condition formulated by a linear
two-dimensional random equation. This equation allows us to
calculate a priori estimate for the fixed point theorem. We note
that this linear random differential equation has a nontrivial
invariant manifold if and only if the gap condition is satisfied.
Hence, our results are
optimal in this sense.\\

We believe  that our technique can be applied to other cases that
are treated in  Bates et al. \cite{BatLuZen98}.\\

We also note that we do not need to use the semigroup given by the skew product flow.\\

  In Section \ref{s2}, we recall some  basic concepts for
random dynamical systems and show that the stochastic partial
differential equation (\ref{spde0})  generates a random dynamical
system. We introduce a random graph transform in Section \ref{s3}.
A generalized   fixed point theorem is presented in Section
\ref{ss}. Finally, we present the main theorem  on invariant manifolds
 in Section \ref{s4}.

%%%%%%%%%%%%%%%%%%%%%
%%%%%%%%%%%%%%%%%%%%%
%%%%  Section 2     %%%%%%%%%%
%%%%%%%%%%%%%%%%%%%%%
\section{Random dynamical systems}\label{s2}

 We recall some basic concepts in
  random dynamical systems. Let
$(\Omega,\mathcal{F},\mathbb{P})$ be a probability space.
A {\em flow} $\theta$ of mappings
$\{\theta_t\}_{t\in\mathbb{R}}$ is defined  on the sample
space  $\Omega$   such that
\begin{equation}\label{eqc1}
\theta:\mathbb{R}\times \Omega\to \Omega,\qquad \theta_0={\rm
id}_\Omega, \qquad\theta_{t_1}\circ\theta_{t_2}= \theta_{t_1+t_2}
\end{equation}
for $t_1,\,t_2\in\mathbb{R}$. This flow is supposed to be
$(\mathcal{B}(\mathbb{R})\otimes\mathcal{F},\mathcal{F})$-measurable,
where $\mathcal{B}(\mathbb{R})$ is the collection of Borel sets on the real line $ \mathbb{R}$.
To have this measurability,  it is not allowed to replace
$\mathcal{F}$ by its $\mathbb{P}$-completion
$\mathcal{F}^\mathbb{P}$; see Arnold \cite{Arn98} Page 547. In
addition, the measure $\mathbb{P}$ is assumed to be ergodic with
respect to $\{\theta_t\}_{t\in\mathbb{R}}$.
Then $(\Omega,\mathcal{F},\mathbb{P},\mathbb{R},\theta)$ is called a metric dynamical system.\\

For our applications, we will consider a special but very important metric
dynamical system  induced by the {\em Brownian motion}. Let $W(t)$ be
a two-sided Wiener process with trajectories in the space
$C_0(\mathbb{R},\mathbb{R})$ of real
continuous functions defined on $\mathbb{R}$,
  taking zero value  at $t=0$. This set is
equipped with the compact open topology. On this set we consider
the measurable flow
$\theta=\{\theta_t\}_{t\in\mathbb{R}},\, \mbox{defined by} \,\theta_t\omega=\omega(\cdot+t)-\omega(t)$.
The distribution of this process generates a measure on
$\mathcal{B}(C_0(\mathbb{R},\mathbb{R}))$ which is called the {\em
Wiener measure}. Note that this measure is ergodic with respect to
the above flow; see the Appendix in Arnold \cite{Arn98}. Later on
we will consider, instead of the whole $C_0(\mathbb{R},\mathbb{R})$,
 a
$\{\theta_t\}_{t\in\mathbb{R}}$-invariant subset $\Omega\subset
C_0(\mathbb{R},\mathbb{R})$ of $\mathbb{P}$-measure one and the
trace $\sigma$-algebra $\mathcal{F}$ of
$\mathcal{B}(C_0(\mathbb{R},\mathbb{R}))$ with respect to
$\Omega$. A set $\Omega$ is called $\{\theta_t\}_{t\in\mathbb{R}}$-invariant
if $\theta_t\Omega=\Omega$ for $t\in\mathbb{R}$.
On $\mathcal{F}$ we consider the restriction of the
Wiener measure also denoted by $\mathbb{P}$.\\

The dynamics of the system on the state space $H$ over the flow
$\theta$ is described by a cocycle. For our applications it is
sufficient to assume that $(H,d_H)$ is a complete metric space. A
cocycle $\phi$ is a mapping:
\[
\phi:\mathbb{R}^+\times \Omega\times H\to H
\]
which  is
$(\mathcal{B}(\mathbb{R})\otimes\mathcal{F}\otimes\mathcal{B}(H),\mathcal{F})$-measurable
 such that
\begin{equation*}
%\label{eq3}
\begin{split}
&\phi(0,\omega,x)=x \in H,\\
&
\phi(t_1+t_2,\omega,x)=\phi(t_2,\theta_{t_1}\omega,\phi(t_1,\omega,x)),
\end{split}
\end{equation*}
for $t_1,\,t_2\in\mathbb{R}^+,\,\omega\in \Omega,$ and $x\in H$.
Then $\phi$ together with the metric dynamical system forms a {\em
random dynamical system}.
\\

   Random dynamical systems are usually generated by
differential equations with random coefficients
\[
\phi^\prime=f(\theta_t\omega,\phi),\quad \phi(0)=x\in\mathbb{R}^d
\]
or finite dimensional stochastic differential equations
\[
d\phi=f(\phi)dt+g(\phi)dW,\quad \phi(0)=x\in\mathbb{R}^d
\]
provided that the global existence and the uniqueness can be
ensured. For details see Arnold \cite{Arn98}. We call a random
dynamical system {\em continuous} if the mapping
\[
x\to\phi(t,\omega,x)
\]
is continuous for $t\in\mathbb{R}^+$ and $\omega\in\Omega$.\\

Now we  start our  investigation on the following stochastic
partial differential equation
\begin{equation}\label{u3}
\frac{d\phi}{dt}=A\phi+F(\phi)+\phi\,\dot{W}
\end{equation}
on a separable Banach space $(H,\|\cdot\|_H)$. Here $A$ is a
linear partial differential operator;  $W(t)$ is an one
dimensional standard Wiener process, and $\dot{W}$ describes
formally a {\em white noise}.
 Note that $\phi\,\dot{W}$ is interpreted
as a Stratonovich differential. However, the existence theory for
stochastic evolution equations is usually formulated for Ito
equations as in Da Prato and Zabczyk \cite{DaPZab92}, Chapter 7.
The equivalent Ito equation for (\ref{u3}) is given by
\begin{equation}\label{eq3a}
d\phi=A\phi\,dt+F(\phi)\,dt+\frac{\phi}{2}dt+\phi\,dW.
\end{equation}

In the following, we assume that the linear (unbounded) operator
$A: D(A)\to H$ generates a strongly continuous semigroup
$\{S(t)\}_{t\ge 0}$ on $H$. Furthermore, we assume that $S(t)$
satisfies the exponential dichotomy with exponents $\hat \lambda
> \check \lambda$ and bound $M$, i.e., there exists a
continuous projection $\pi^+$ on $H$ such that
\begin{itemize}
\item[(i)] $\pi^+S(t)=S(t)\pi^+$;
\item[(ii)] the restriction $S(t)|_{R(\pi^+)}$, $t\geq 0$, is an
isomorphism of ${R(\pi^+)}$ onto itself, and we define $S(t)$ for
$t <0$ as the inverse map.
\item[(iii)]
\begin{equation}\label{eqvv}
\|\pi^+S(t)\pi^+\|_{H,H}\le Me^{\hat \lambda t},\quad t\le 0,\quad
\|\pi^-S(t)\pi^-\|_{H,H}\le Me^{\check \lambda t},\quad t\ge 0\;
\end{equation}
where $\pi^-=I-\pi^+$.
\end{itemize}
Denote $H^-= \pi^-H$ and $H^+=\pi^+H$. Then, $H=H^+\oplus H^-$.
\\

For simplicity we set  $M=1$. For instance, if the operator $-A$ is a
strongly elliptic and symmetric differential operator on a smooth
domain  $\bar D$ of order $2$   under the homogeneous
Dirichlet boundary conditions, then the above assumptions are
satisfied with $H=L^2(D)$. In this case $A$ has the spectrum
\begin{equation*}
%\label{eq441}
\lambda_1>\cdots>\lambda_u>\lambda_{u+1} > \lambda_{u+2}> \cdots
\end{equation*}
where the space spanned by the associated eigenvectors is equal to
$H$. For any $\lambda_u$ the associated eigenspace is finite
dimensional. The space $H^+$  is spanned by the associated
eigenvectors for $\lambda_{1},\,\lambda_{2},\cdots,\lambda_u$
and $\hat\lambda=\lambda_u>\lambda_{u+1}=\check\lambda$.\\

We assume that  $F$ is  Lipschitz continuous on $H$
\[
\|\pi^\pm (F(x_1)-F(x_2))\|_H\le L\|x_1-x_2\|_H
\]
with the Lipschitz constant $L>0$. Then, for any initial data
$x\in H$, there exists a unique solution of (\ref{eq3a}). For
details about the properties of this solution see Da
Prato and Zabczyk \cite{DaPZab92}, Chapter 7.
We also assume that $F(0)=0$.\\

The   stochastic evolution equation (\ref{eq3a}) can be written   in the
 following {\em mild}  integral form:
\[
\phi(t)=S(t)x+\int_0^t(S(t-\tau)(F(\phi(\tau))+\frac{\phi(\tau)}{2})d\tau
+\int_0^tS(t-\tau)\phi(\tau)dW,\quad x\in H
\]
almost surely for any $x\in H$. Note that the theory in
\cite{DaPZab92}
requires that the associated probability space is  complete.\\

In order to apply the random dynamical systems techniques, we
introduce a coordinate transform converting conjugately a
stochastic  partial differential equation into an infinite
dimensional random dynamical system. Although it is well-known
that a large class of partial differential equations with
stationary random coefficients and  Ito stochastic ordinary
differential equations generate random dynamical systems (for
details see Arnold \cite{Arn98}, Chapter 1), this problem is still
unsolved for stochastic partial differential equations with a
general noise term $C(u)\,dW$. The reasons are: (i) The stochastic
integral is only defined almost surely where the exceptional set
may depend on the initial state $x$;   and (ii)
  Kolmogorov's theorem,   as cited in
Kunita \cite{Kun90} Theorem 1.4.1,  is only true for finite
dimensional random fields.
Moreover,  the cocycle has to be defined for  {\em any} $\omega\in\Omega$.\\
However, for the noise term $\phi\,dW$ considered here,  we can
show that (\ref{eq3a}) generates a random dynamical system. To
prove this property, we need the following preparation.\\

We consider the one-dimensional linear stochastic differential equation:
\begin{equation}\label{equ4}
dz+z\,dt=dW.
\end{equation}
A solution of this equation is called an Ornstein-Uhlenbeck
process.
\begin{lemma}\label{l100}
i) There exists a $\{\theta_t\}_{t\in\mathbb{R}}$-invariant set
$\Omega\in\mathcal{B}(C_0(\mathbb{R},\mathbb{R}))$ of full measure
with sublinear growth:
\[
\lim_{t\to\pm\infty}\frac{|\omega(t)|}{|t|}=0,\qquad
\omega\in\Omega
\]
of $\mathbb{P}$-measure one.\\
ii) For $\omega\in\Omega$ the random variable
\[
z(\omega)=-\int_{-\infty}^0e^\tau\omega(\tau)d\tau
\]
exists and generates a unique stationary solution of (\ref{equ4})
given by
\[
\Omega\times \mathbb{R}\ni(\omega,t)\to z(\theta_t\omega)
=-\int_{-\infty}^0e^\tau\theta_t\omega(\tau)d\tau
=-\int_{-\infty}^0e^\tau\omega(\tau+t)d\tau+\omega(t).
\]
The mapping $t\to z(\theta_t\omega)$
is continuous.\\
iii) In particular, we have
\begin{equation*}
%\label{eq401}
\lim_{t\to\pm\infty}\frac{|z(\theta_t\omega)|}{|t|}=0\quad
\text{for }\omega\in\Omega.
\end{equation*}
iv) In addition,
\[
\lim_{t\to\pm\infty}\frac{1}{t}\int_0^tz(\theta_\tau\omega)d\tau=0
\]
for $\omega\in\Omega$.
\end{lemma}
\begin{proof}
i) It follows from the law of iterated logarithm that there exists
a set
$\Omega_1\in\mathcal{B}(C_0(\mathbb{R},\mathbb{R})),\,\mathbb{P}(\Omega_1)=1$,
such that
\[
\limsup_{t\to\pm\infty}\frac{|\omega(t)|}{\sqrt{2|t|\log\log
|t|}}=1
\]
for  $\omega\in\Omega_1$.
The set of these $\omega$'s is $\{\theta_t\}_{t\in\mathbb{R}}$-invariant.\\

ii) This can be proven as in {\O}ksendal \cite{Oks92} Page 35. The
existence of the integral on the right hand side for
$\omega\in\Omega_1$ follows from the law of iterated logarithm.
Using the law of iterated logarithm again, the function
\[
\tau\to e^\tau\sup_{[t_0-1,t_0+1]}|\omega(\tau+t_0)|
\]
is an integrable majorant for $e^\tau \omega(\tau+t)$ for $t\in
[t_0-1,t_0+1]$ and $\tau\in (-\infty,0]$. Hence the continuity at
$t_0\in\mathbb{R}$
follows straightforwardly from Lebesgue's theorem of dominated convergence.\\

iii) By the law of iterated logarithm,  for $1/2<\delta<1$ and
$\omega\in\Omega_1$ there exists a constant $C_{\delta,\omega}>0$
such that
\[
|\omega(\tau+t)| \le C_{\delta,\omega}+|\tau+t|^\delta\le
C_{\delta,\omega}+|\tau|^\delta+|t|^\delta,\quad \tau\le 0.
\]
Hence
\begin{align*}
&\lim_{t\to\pm\infty}
\left|\frac{1}{t}\int_{-\infty}^0e^\tau\omega(\tau+t)d\tau\right|\le
\lim_{t\to\pm\infty}
\frac{1}{|t|}\int_{-\infty}^0e^\tau (C_{\delta,\omega}+|\tau|^\delta+|t|^\delta)d\tau=0,\\
&\lim_{t\to\pm\infty}\frac{\omega(t)}{t}=0
\end{align*}
which gives  the convergence relation in iii). Hence, these
convergence relations  always define a
$\{\theta_t\}_{t\in\mathbb{R}}$-invariant set which has a full measure.\\

iv) Clearly, $\mathbb{E}z=0$ from ii). Hence by the ergodic
theorem we obtain iv) for
$\omega\in\Omega_2\in\mathcal{B}(C_0(\mathbb{R},\mathbb{R}))$.
This set  $\Omega_2$ is also
$\{\theta_t\}_{t\in\mathbb{R}}$-invariant. Then we set
\[
\Omega:=\Omega_1\cap\Omega_2.
\]
The proof is complete.
\end{proof}

We now replace $\mathcal{B}(C_0(\mathbb{R},\mathbb{R}))$ by
\[
\mathcal{F}=\{\Omega\cap A,\,A\in
\mathcal{B}(C_0(\mathbb{R},\mathbb{R}))\}
\]
for $\Omega$ given in Lemma \ref{l100}. The probability measure is
the restriction of the Wiener measure to this new
$\sigma$-algebra, which is  also denoted by $\mathbb{P}$. In the
following we will consider the metric dynamical system
\[
(\Omega,\mathcal{F},\mathbb{P},\mathbb{R},\theta).
\]

We now back to show that the solution of (\ref{eq3a}) defines a random
dynamical system. To see this, we
 consider the random partial differential equation
\begin{equation}\label{u5}
\frac{d\phi}{dt}=A\phi+G(\theta_t\omega,\phi)+z(\theta_t\omega)\phi,\quad
\phi(0)=x\in H
\end{equation}
where $G(\omega,u):=e^{-z(\omega)}F(e^{z(\omega)}u)$. It is easy
to see that for any $\omega\in\Omega$ the function $G$ has the
same global Lipschitz constant $L$ as $F$. In contrast to the original
stochastic differential equation, no stochastic integral appears here.
The solution can be interpreted in a mild sense
\begin{equation}\label{u5a}
\phi(t)=e^{\int_0^tz(\theta_\tau\omega)d\tau}S(t)x+\int_0^t
e^{\int_\tau^tz(\theta_r\omega)dr}
S(t-\tau)G(\theta_\tau\omega,\phi(\tau))d\tau.
\end{equation}
We note that this equation has  a unique solution for every
$\omega\in\Omega$. No exceptional sets appear.
 Hence the solution mapping
\[
(t,\omega,x)\to\phi(t,\omega,x)
\]
generates a random dynamical system. Indeed, the mapping $\phi$ is
$(\mathcal{B}(\mathbb{R})\otimes\mathcal{F}\otimes\mathcal{B}(H),\mathcal{F})$-measurable.\\

Let   $\hat\phi(t,\omega,x)$ be the solution mapping of
(\ref{eq3a}) which is defined for
$\omega\in\Omega\in\mathcal{F}^\mathbb{P},\,\mathbb{P}(\Omega)=1$.
We now introduce the transform
\begin{equation}\label{eq500}
T(\omega,x)=xe^{-z(\omega)}
\end{equation}
and its inverse transform
\begin{equation}\label{eq501}
T^{-1}(\omega,x)=xe^{z(\omega)}
\end{equation}
for $x\in H$ and $\omega\in\Omega$.

\begin{lemma}\label{l111}
Suppose that $\phi$ is the random dynamical system generated by
(\ref{u5}).     Then
\begin{equation}\label{eq503}
(t,\omega,x)\to T^{-1}(\theta_t\omega,\cdot)\circ
\phi(t,\omega,T(\omega,x))=:\hat\phi(t,\omega,x)
\end{equation}
is a random dynamical system. For any $x\in H$ this process is a
solution version of (\ref{eq3a}).
\end{lemma}

\begin{proof}
Applying the Ito formula to
$T(\theta_t\omega,\hat\phi(t,\omega,T^{-1}(\omega,x)))$  gives a
solution of (\ref{u5}). The converse is also true, since
$T^{-1}(\theta_t\omega,x)$ and $\phi(t,\omega,x))$ are defined for
{\em any} $\omega\in\Omega$ and $T^{-1}$ is the inverse  of $T$, and thus
\[
(t,\omega,x)\to T^{-1}(\theta_t\omega,\phi(t,\omega,T(\omega,x)))
\]
gives a solution of (\ref{eq3a})  for each $\omega\in\Omega$. It
is easy to check that (\ref{eq503}) defines a random dynamical
system. Since $\phi$ is measurable with respect to $\mathcal{F}$
so is this $\hat\phi$.
\end{proof}

Similar transformations     have been used by Caraballo, Langa and
Robinson \cite{CarLanRob01} and  Schmalfu{\ss} \cite{Schm97c}. Note
that our transform has the  advantage that the solution of
(\ref{u5}) generates a random dynamical system for the
$\omega$-wise differential equation.\\

In Section \ref{s4} we will prove the existence of  invariant
manifolds generated by (\ref{u5}). These  invariant  manifolds can be
transformed into invariant manifolds for (\ref{u3}).

%%%%%%%%%%%%%%%%%
%%%%%%%%%%%%%%%%%
%%%%%%%%%%%%%%%%%
\section{Random graph transform}\label{s3}

In this section, we construct a random graph transform. The fixed
point of this transform gives the desired invariant manifold for
the random dynamical system $\phi$ generated by (\ref{u5}).\\

We first recall that a multifunction
$M=\{M(\omega)\}_{\omega\in\Omega}$ of nonempty closed sets
$M(\omega),\,\omega\in\Omega$, contained in a complete separable
metric space $(H,d_H)$ is called a {\em random set} if
\[
\omega\to\inf_{y\in M(\omega)}d_H(x,y)
\]
is a random variable for any $x\in H$.
\begin{definition}
A random set $M(\omega)$ is called an invariant set if
\[
\phi(t,\omega,M(\omega))\subset M(\theta_t\omega).
\]
If we can represent $M$ by a graph of a Lipschitz mapping
\[
\gamma^\ast(\omega,\cdot): \;  H^+\to H^-,\quad H^+\oplus H^-=H
\]
such that
\begin{equation*}
%\label{eq1200}
M(\omega)=\{x^++\gamma^\ast(\omega,x^+),x^+\in H^+\},
\end{equation*}
then $M(\omega)$ is called a Lipschitz continuous invariant
manifold.
\end{definition}

Let $\gamma(\cdot):H^+\to H^-$ be a  Lip\-schitz continuous
function with Lipschitz constant $L_\gamma\ge 0$ and also  let
$\gamma(0)=0$. We consider the system of equations
\begin{align}\label{eq8}
\begin{split}
w(t)=&e^{\int_T^tz(\theta_\tau\omega)d\tau}\pi^+S(t-T)y^+\\
&-\int_t^T
e^{\int_\tau^tz(\theta_{\tau^\prime}\omega)d\tau^\prime}\pi^+S(t-\tau)\pi^+G(\theta_\tau\omega,w(\tau)+v(\tau))d\tau\\
v(t)=&e^{\int_0^tz(\theta_\tau\omega)d\tau}\pi^-S(t)\gamma(w(0))\\
&+\int_0^te^{\int_\tau^tz(\theta_{\tau^\prime}\omega)d\tau^\prime}\pi^-S(t-\tau)\pi^-G(\theta_\tau\omega,w(\tau)+v(\tau))d\tau
\end{split}
\end{align}
on some interval $[0,T]$. Note that if (\ref{eq8}) has a solution
$(w,v)$ on $[0,T]$ then $w(0)$ defines a mapping
$\gamma\to\Psi(T,\theta_T\omega,\gamma)(y^+)$ and $v(T)$
defines another mapping
\begin{equation}\label{eq7}
\gamma\to\Phi(T,\omega,\gamma)(y^+).
\end{equation}
This latter mapping  $\Phi$  will serve as the random graph transform.
\\

Recall that a random variable $\omega\to\gamma^\ast(\omega)$ is a generalized fixed point
of the mapping $\Phi$ if
\begin{equation}
\Phi(T,\omega,\gamma^\ast(\omega))=\gamma^\ast(\theta_T\omega).
\end{equation}
for $\omega\in\Omega,\,T\ge 0$. We assume that
$\gamma^\ast(\omega)$ is a Lipschitz continuous mapping
from $H^+$ to $H^-$ and it
  takes zero  value at zero. Conditions for the existence of a generalized
fixed point are derived in the next section in the case of
$\Phi$  a random dynamical system. The following theorem describes
the relation between generalized fixed points and invariant manifolds.
\begin{theorem}\label{main}
Suppose that $\gamma^\ast$ is the generalized fixed point of the
mapping $\Phi$.
Then the graph of   $\gamma^\ast$  is the  invariant manifold $M(\omega)$
of the random dynamical system $\phi$
generated by (\ref{u5}).
\end{theorem}
\begin{proof}
Let $M(\omega)$ be the graph of $\gamma^\ast(\omega)$ such that $(x^+,\gamma^\ast(x^+,\omega))\in M(\omega)$.
Then     for $x^+,\,y^+\in H^+$, we obtain
\begin{align*}
\phi(T,\omega,x^++\gamma^\ast(\omega,x^+))=&
\pi^+\phi(T,\omega,x^++\gamma^\ast(\omega,x^+))+
\pi^-\phi(T,\omega,x^++\gamma^\ast(\omega,x^+))\\
=&
y^++\pi^-\phi(T,\omega,\Psi(T,\theta_T\omega,\gamma^\ast(\omega))(y^+)\\
&
+\gamma^\ast(\omega,\Psi(T,\theta_T\omega,\gamma^\ast(\omega))(y^+)))\\
=&y^++\Phi(T,\omega,\gamma^\ast(\omega))(y^+)=y^++\gamma^\ast(\theta_T\omega)(y^+)\in
M(\theta_T\omega)
\end{align*}
by the definition of $\Psi$:
\[
x^+=\Psi(T,\theta_T\omega,\gamma^\ast(\omega))(y^+)\;\text{if and only if }
y^+=\pi^+\phi(T,\omega,x^++\gamma^\ast(\omega,x^+)).
\]
For the measurability statement see Section \ref{s4} below.
\end{proof}

By this theorem, we can find invariant  manifolds of the
random dynamical system $\phi$ generated by  (\ref{u5})
by finding generalized fixed points of the mapping  $\Phi$ defined
in (\ref{eq7}). To do so, we will use a generalized
fixed point theorem for {\em cocycles} and thus we need
to show that the above mapping  $\Phi$ is in fact a random
dynamical system.
For the remainder of this section we will show that $\Phi$ defines a random dynamical system. We will achieve this in a few lemmas.\\

In the following we denote by $C^{0,1}_0(H^+;B)$ the Banach space
of Lipschitz continuous functions from $H^+$, with value zero at
zero, into a Banach space $B$ with the usual   (Lipschitz)  norm
\[
\|u\|_{C_0^{0,1}}=\sup_{y_1^+\not= y_2^+\in
H^+}\frac{\|u(y_1^+)-u(y_2^+)\|_B}{\|y_1^+-y_2^+\|_H}.
\]
Moreover, $C_0^{G}(H^+;B)$ denotes the Banach space of
bounded continuous functions, with value zero at zero and
with linearly growth. The norm in this space is defined as
\[
\|u\|_{C_0^G}=\sup_{0\not=y^+\in
H^+}\frac{\|u(y^+)\|_B}{\|y^+\|_H}.
\]

 We first present a result about
the existence of a solution of the integral system (\ref{eq8}).
 The proof is quite technical and is given in the Appendix.\\

\begin{lemma}\label{l17}
Let $L$ be the Lipschitz constant of the nonlinear
term $G(\omega,\cdot)$   in the random partial differential equation
(\ref{u5}). Then for any $\gamma\in
C_0^{0,1}(H^+;H^-),\,\omega\in\Omega$,    there exists  a
$T=T(\gamma,\omega)>0$ such that on $[0,T]$   the integral system
(\ref{eq8})   has a unique
solution    $(w(\cdot),v(\cdot))\in
C([0,T];C_0^{G}(H^+;H^+)\times C_0^{G}(H^+;H^-))$.
\end{lemma}

\medskip

Let $C([0,T];B)$ be the space of continuous mappings from $[0,T]$ into
$B$.
Note that for some $T>0$ and $\gamma\in C_0^{0,1}(H^+;H^-)$, the fixed
point problem defined by the integral system
(\ref{eq8}) has a contraction constant less than one. Then for
$T^\prime<T$ and some Lipschitz continuous function
$\gamma^\prime\in C_0^{0,1}(H^+;H^-)$ such that
$\|\gamma^\prime\|_{C_0^{0,1}}\le \|\gamma\|_{C_0^{0,1}}$ the same
contraction constant can be chosen. This follows from the
structure of the contraction constant; see (\ref{000}) below.
\\

We would like to calculate a priori estimates for
the solution of (\ref{eq8}). To do this we   need the following
lemma and its conclusion on monotonicity will  also be used later on.

\begin{lemma}\label{l18}
We consider the differential equations
\begin{align}\label{eq11}
\begin{split}
W^\prime &=\hat\lambda W+z(\theta_t\omega)W-LW-LV,\qquad \\
V^\prime& =\check\lambda V+z(\theta_t\omega)V+LW+LV
\end{split}
\end{align}
with generalized initial conditions
\begin{align}\label{eq11b}
W(T)=Y\ge 0,\qquad   V(0)=\Gamma W(0)+C,\;\Gamma,\,C\ge 0.
\end{align}
Then this system has a unique solution  on $[0,T]$ for  some
$T=T(\Gamma,C,\omega)>0$. This interval is independent of $C$. Let
$\hat W,\,\hat V$ be  solutions of (\ref{eq11}) but with the
generalized initial conditions
\[
\hat W(T)=Y\ge 0, \quad \hat V(0)=\hat\Gamma\hat W(0)+\hat
C,\qquad 0\le \hat\Gamma\le\Gamma,\;0\le \hat C\le C.
\]
Then we have $0\le \hat V(t)\le V(t)$ and $0\le \hat W(t)\le W(t)$ for
$t\in [0,T]$.
\end{lemma}
The proof is given in the Appendix.\\

Now we can compare the norms for the solution of (\ref{eq8}) and that of
  (\ref{eq11})-(\ref{eq11b}).

\begin{lemma}\label{l55}
Let $[0,T]$ be an interval on which the assumptions of the Banach
fixed point theorem (see the proofs of Lemma \ref{l17}, \ref{l18})
are satisfied for (\ref{eq8}) and (\ref{eq11})-(\ref{eq11b}) for
some $\gamma\in C_0^{0,1}(H^+;H^-)$. Then the norm of the solution
of (\ref{eq8}) is bounded by the solution of (\ref{eq11})-
(\ref{eq11b}) with $Y=1$,  $C=0$, and $\Gamma=L_\gamma$  being
the Lipschitz norm of $\gamma$. That is,
\[
\|w(t)\|_{C_0^G}\le W(t),\quad \|v(t)\|_{C_0^G}\le V(t).
\]
\end{lemma}

The proof is given in the Appendix.\\

We obtain from Lemma \ref{l17} that $w(t,y^+),\,v(t,y^+)$ exist
for any $y^+\in H^+$ on some interval $[0,T]$. We also   have $\|w(T)\|_{C_0^{0,1}}=1$ and
\begin{align*}
\frac{\|\gamma(w(0,y_1^+))-\gamma(w(0,y_2^+))\|_H}{\|y_1^+-y_2^+\|_H}=&
\frac{\|\gamma(w(0,y_1^+))-\gamma(w(0,y_2^+))\|_H}
{\|w(0,y_1^+)-w(0,y_2^+)\|_H}\times\\
&\times\frac{\|w(0,y_1^+)-w(0,y_2^+)\|_H}{\|y_1^+-y_2^+\|_H}\\
\le &L_\gamma\|w(0)\|_{C_0^{0,1}}
\end{align*}
for $y_1^+\not=y_2^+$ and $w(0,y_1^+)\not=w(0,y_2^+)$.  Hence
$\|v(0)\|_{C_0^{0,1}}\le L_\gamma \|w(0)\|_{C_0^{0,1}}$. We have
that $w(0,y_1^+)\not=w(0,y_2^+)$ because
$\Psi(T,\theta_T\omega,\gamma)(\cdot)$ is a bijection. Indeed this
mapping is the inverse of
$x^+\to\pi^+\phi(T,\omega,x^++\gamma(x^+))$ on $H^+$. One can see
this if we plug in $x^+=\Psi(T,\theta_T\omega,\gamma)(\cdot)$,
which is given by $w(0)$, the right hand side of (\ref{eq8}) at
zero into the $\pi^+$-projection of the right hand side of
(\ref{u5a}) for $t=T$, and vice versa if we plug in this expression
into the right hand side of the first equation of (\ref{eq8}). On
the other hand,  we have
\begin{align*}
&\frac{\|\pi^\pm G(\omega,w(y_1^+)+v(y_1^+))-\pi^\pm G(\omega,w(y_2^+)+v(y_2^+))\|_H}{\|y_1^+-y_2^+\|_H}\\
&\quad\le
L\frac{\|w(y_1^+)+v(y_1^+)-(w(y_2^+)+v(y_2^+))\|_H}{\|y_1^+-y_2^+\|_H}\\
&\quad\le L\frac{\|w(y_1^+)-w(y_2^+)\|_H}{\|y_1^+-y_2^+\|_H}
+L\frac{\|v(y_1^+)-v(y_2^+)\|_H}{\|y_1^+-y_2^+\|_H}.
\end{align*}
Repeating the arguments of Lemma \ref{l55} we obtain
\[
\frac{\|w(y_1^+)-w(y_2^+)\|_H}{\|y_1^+-y_2^+\|_H}\le W(t),\quad
\frac{\|v(y_1^+)-v(y_2^+)\|_H}{\|y_1^+-y_2^+\|_H}\le V(t)
\]
for any $y_1^+\not= y_2^+$. Hence, we have the following result.

\begin{lemma}\label{l55a}
The solution of  the integral system
(\ref{eq8})  has the  following regularity:
 $w(t)\in C_0^{0,1}(H^+;H^+)$ and
$v(t)\in C_0^{0,1}(H^+;H^-)$. In particular,
$\Phi(T,\omega,\gamma)\in C_0^{0,1}(H^+;H^-)$ for sufficiently
small $T$. Moreover, the comparison result in Lemma \ref{l55}
remains true.
\end{lemma}

Note that by the fixed point argument,  $\Phi(T,\omega,\gamma)$ and
$\Psi(T,\theta_T\omega,\gamma)$ exist only for small $T$. We would
like to extent these definitions  to $T\in \mathbb{R}^+$. To see
this, we are going to show that if the Lipschitz constant of
$\gamma$ is bounded by a particular value, then the Lipschitz
constant of
 $\mu=\Phi(T,\omega,\gamma)$
has the same bound.\\
As a preparation we consider the matrix
\[
B:=\begin{pmatrix}
\hat\lambda-L&-L\\
L&\check\lambda +L
\end{pmatrix}
\]
which has the  eigenvalues $\lambda_+,\lambda_-$. These
eigenvalues are real and distinct if and only if
\begin{equation}\label{eq30}
\hat\lambda-\check\lambda >4L.
\end{equation}
Then  the associated eigenvectors can be written as
\[
(e_+,1),\,(e_-,1).
\]
 We order $\lambda_+,\,\lambda_-$ as
$\lambda_+>\lambda_-$. The elements $e_+,\,e_-$ are positive.
\begin{lemma}\label{l20}
Let $T=T(\Gamma,0,\omega)>0$ be chosen such that (\ref{eq11}),
(\ref{eq11b}) have a solution on $[0,T]$ given by the fixed point
argument for $\Gamma=e_+^{-1}=:\kappa,\,Y=1$ and $C=0$. Then the
closed ball $B_{C_0^{0,1}}(0,\kappa)$ in $C_0^{0,1}(H^+;H^-)$ will
be mapped into itself:
$\Phi(T,\omega,B_{C_0^{0,1}}(0,\kappa))\subset
B_{C_0^{0,1}}(0,\kappa)$.
\end{lemma}
\begin{proof}
Let $Q_1(t)\vec{x}_0$ be the solution  of the linear initial value
problem
\[
\vec{x}^\prime=B\vec{x},\quad\vec{x}(0)=\vec{x}_0
\]
and let
\[
Q_2(t)=
\begin{pmatrix}
e^{\int_0^tz(\theta_\tau\omega)d\tau}&0\\
0&e^{\int_0^tz(\theta_\tau\omega)d\tau}
\end{pmatrix}
\]
be the solution operator of
\begin{align*}
\psi^\prime&=z(\theta_t\omega)\psi,\quad\psi(0)=\psi_0,\\
\eta^\prime&=z(\theta_t\omega)\eta,\quad\eta(0)=\eta_0.
\end{align*}
Note that $Q_2(t)$ and $Q_1(t)$ commute. Hence $Q_2(t)Q_1(t)$ is a
solution operator of the linear differential equation
(\ref{eq11}). Since
\[
Q_1(t) \left(
\begin{array}{c}
e_+\\1
\end{array}
\right)=e^{\lambda_+ t}\left(
\begin{array}{c}
e_+\\1
\end{array}
\right)
\]
we obtain that
\[
Q_2(t)Q_1(t) \left(
\begin{array}{c}
e_+\\1
\end{array}
\right)=e^{\lambda_+ t+\int_0^tz(\theta_\tau\omega)d\tau}\left(
\begin{array}{c}
e_+\\1
\end{array}
\right).
\]
For the initial conditions $Y=1,\,\Gamma=e_+^{-1}$ we can calculate
explicitly for the solution of (\ref{eq11}), (\ref{eq11b})
\[
W(0)=e^{-\lambda_+T-\int_0^Tz(\theta_\tau\omega)d\tau},\quad
c_1=e_+^{-1}e^{-\lambda_+T-\int_0^Tz(\theta_\tau\omega)d\tau}
\]
and $c_2=0$. Hence  $V(T)=e_+^{-1}$. By  the
comparison results from Lemmas \ref{l55} and\ref{l55a}, we find that
$\|w(0)\|_{C_0^{0,1}}\le W(0)$ and
$\|v(T)\|_{C_0^{0,1}}=\|\Phi(T,\omega,\gamma)\|_{C_0^{0,1}}\le
V(T)=e_+^{-1}$ for small $T$ depending on $\omega$ such that
\[
\Phi(T,\omega,B_{C_0^{0,1}}(0,\kappa))\subset
B_{C_0^{0,1}}(0,\kappa).
\]
\end{proof}

 Since we will equip $B_{C_0^{0,1}}(0,\kappa)$
with the $C_0^G$-norm in Section  \ref{s4},  in the following we will choose the state
space
$\mathcal{H}=B_{C_0^{0,1}}(0,\kappa)$ with the metric  $d_\mathcal{H}(x,y):=\|x-y\|_{C_0^G}$.\\

Now we show that the random graph transform $\Phi$ defines a
random dynamical system.

\begin{theorem}\label{t11}
Suppose that the gap condition (\ref{eq30}) is satisfied. Then
$\Phi$ is  well-defined by (\ref{eq7}) for {\em any} $T\ge
0,\,\omega\in\Omega$ and $\gamma\in\mathcal{H}$. In addition,
$\Phi$ together with the metric dynamical system $\theta$ induced
 by the
Brownian motion defines a random dynamical system. In particular,
the following measurability for the operators of the cocycle
holds:
\[
\Omega\ni\omega\to \Phi(T,\omega,\gamma)(y^+)\in H^-
\]
is $(\mathcal{F},\mathcal{B}(H^-))$-measurable for any $y^+\in
H^+,\,T\ge 0$.
\end{theorem}
\begin{proof}
By Lemma \ref{l17}, the mapping $\Phi(T,\omega,\gamma)$ is  defined for
small $T$. So we first   have to extend this definition for any
$T>0$.
\\

To this end we introduce random variables $T_\kappa(\omega)>0$ by
\[
T_\kappa(\omega):=\frac{1}{2}\inf\{T>0:K(\omega,T,\kappa)\ge 1\}
\]
where $K$ is defined in (\ref{000})  below. Since $T\to K(\omega,T,\kappa)$ is continuous in $T$ this is a random variable.
Hence,
$K(\omega,T_\kappa(\omega),\kappa)<1$, and  (\ref{eq8}) has a
unique solution on $[0,T_\kappa(\omega)]$ for
$\gamma\in\mathcal{H}$. We define a sequence by
$T_1=T_1(\omega)=T_\kappa(\omega)$,
$T_2=T_2(\omega)=T_\kappa(\theta_{T_1(\omega)}\omega)$ and so on.
Suppose that for some $\omega\in\Omega$ we have that
$\sum_{i=1}^\infty T_i(\omega)=T_0<\infty$. Then  the definition
of $K$ in (\ref{000}) implies that
$\int_0^{T_0}|z(\theta_\tau\omega)|d\tau=\infty$. This
 is a contradiction, because
by Lemma \ref{l100} the mapping $t\to z(\theta_t\omega)$ is
continuous. Hence for any $T>0$ and $\omega\in\Omega$ there exists
an $i=i(T,\omega)$ such that
\[
T=T_1+T_2+\cdots +T_{i-1}+\hat T_i,\quad 0<\hat T_i\le T_i.
\]
We   can now define
\begin{equation}\label{001}
\Phi(T,\omega,\gamma)=\Phi(\hat T_i,\theta_{T_{i-1}}\omega,\cdot)\circ
\Phi(T_{i-1},\theta_{T_{i-2}}\omega,\cdot)
\circ\cdots\circ\Phi(T_1,\omega,\gamma).
\end{equation}
We show that the right hand side satisfies (\ref{eq8}). \\
Suppose that
$(w^1,v^1)=(w^1(t,\omega,\gamma,y^+),v^1(t,\omega,\gamma,y^+))$ is
given by (\ref{eq8}) on some interval $[0,t_1],\,t_1\le T_1$ for
$\gamma\in\mathcal{H}$. We have
\[
\mu(\cdot):=v^1(t_1,\omega,\gamma,\cdot)=\Phi(t_1,\omega,\gamma)(\cdot)\in\mathcal{H},
\]
see Lemma \ref{l20}.
Similarly,
\[
(w^2,v^2)=(w^2(t,\theta_{t_1}\omega,\mu,z^+),v^2(t,\theta_{t_1}\omega,\mu,z^+))
\]
is given by (\ref{eq8}) on some interval $[0,t_2],\,t_2\le T_2$.
We set
\[
w(t,\omega,\gamma,z^+)=\left\{
\begin{array}{ccc}
w^1(t,\omega,\gamma,w^2(0,\theta_{t_1}\omega,\mu,z^+))&:&t\in [0,t_1]\\
w^2(t-t_1,\theta_{t_1}\omega,\mu,z^+)&:&t\in(t_1,t_1+t_2]
\end{array}
\right..
\]
By the variation of
constants formula on $w$ we have for $t\in[0,t_1]$
\begin{align}\label{bubu}
\begin{split}
e&^{\int_{t_1}^tz(\theta_\tau\omega)d\tau}\pi^+
S(t-t_1)e^{\int_{t_2}^0z(\theta_{\tau+t_1}\omega)d\tau}\pi^+S(-t_2)z^+\\
&-\pi^+S(t-t_1)e^{\int_{t_1}^tz(\theta_\tau\omega)d\tau}
\int_{0}^{t_2}\pi^+S(-\tau)e^{\int_{\tau}^0z(\theta_{r+t_1}\omega)dr}
\pi^+G(\theta_{\tau+t_1}\omega,w^2+v^2)d\tau\\
&
-\int_{t}^{t_1}\pi^+S(t-\tau)e^{\int_{\tau}^tz(\theta_r\omega)dr}\pi^+
G(\omega,w^1+v^1)d\tau\\
=& e^{\int_{t_1+t_2}^tz(\theta_\tau\omega)d\tau}\pi^+S(t-t_1-t_2)
z^+\\
&-\int^{t_1+t_2}_t\pi^+S(t-\tau)e^{\int_{\tau}^tz(\theta_r\omega)dr}\pi^+G(\theta_{\tau}\omega,w+v)d\tau=w(t).
\end{split}
\end{align}
Now we consider the
second equation of (\ref{eq8}) with initial condition
\[
\gamma(w(0))=\gamma(w^1(0,\omega,\gamma,w^2(0,\theta_{t_1}\omega,\mu,z^+))).
\]
Then at $t_1$ we have for the solution of the second equation
\[
v^1(t_1,\omega,\gamma,w^2(0,\theta_{t_1}\omega,\mu,z^+))=\mu(w^2(0,\theta_{t_1}\omega,\mu,z^+))
\]
which is equal to $v^2(0,\theta_{t_1}\omega,\mu,z^+)$. Hence  for
\[
v(t,\omega,\gamma,z^+)=\left\{
\begin{array}{ccc}
v^1(t,\omega,\gamma,w^2(0,\theta_{t_1}\omega,\mu,z^+))&: &t\in [0,t_1]\\
v^2(t-t_1,\theta_{t_1}\omega,\mu,z^+)&: &t\in(t_1,t_1+t_2]
\end{array}
\right.
\]
we can find
\begin{align*}
v(t_1+t_2)=&e^{\int_0^{t_1+t_2}z(\theta_\tau\omega)d\tau}\pi^-S(t_1+t_2)\gamma(w(0))\\
+&\int_0^{t_1+t_2}e^{\int_\tau^{t_1+t_2}z(\theta_{\tau^\prime}\omega)d\tau^\prime}\pi^-S(t_1+t_2-\tau)\pi^-G(\theta_\tau\omega,w(\tau)+v(\tau))d\tau
\end{align*}
which gives us together with (\ref{bubu}) that $(w,v)$ solves
(\ref{eq8}) on $[0,t_1+t_2]$ and
$v(t_1+t_2)=\Phi(t_1+t_2,\omega,\gamma)(z^+)$. Since $\mu\in\mathcal{H}$ so is
$\Phi(t_1+t_2,\omega,\gamma)(z^+)$ by Lemma (\ref{l20}). The extension of the definition of $\Phi$
is correct since we obtain the same value for different
$t_1\in [0,T_1],\,t_2\in [0,T_1]$ whenever $t_1+t_2={\rm const}$. For this uniqueness we
note that $z\to w(0,\omega,\gamma,z^+)$ given by the above formula
is the inverse of $x^+\to\pi^+\phi(t_1+t_2,\omega,x^++\gamma(x^+))$
which is independent of the choice of $t_1$ and $t_2$. This
implied the independence of $v(t_1+t_2)$ on $t_1+t_2={\rm const}$.
By a special choice of $t_1,\,t_2$ (for instance $t_1=T_1,\,t_2=T_2$
and continuing the above
iteration procedure we get (\ref{001}). By this iteration we also obtain that
$\Phi(T,\omega,\gamma)\in
\mathcal{H}$.

For the measurability,  we note that
\[
\Psi(T\wedge T_\kappa(\omega),\theta_{T\wedge
T_\kappa(\omega)},\gamma)(y^+),\quad \Phi(T\wedge
T_\kappa(\omega),\theta_{T\wedge T_\kappa(\omega)},\gamma)(y^+)
\]
are $\mathcal{F},H^\pm$-measurable because these expressions are
given as an $\omega$-wise limit of the iteration of the Banach
fixed point theorem starting with a measurable expression. On the
other hand,
\[
y^+\to\Psi(T\wedge T_\kappa(\omega),\theta_{T\wedge
T_\kappa(\omega)},\gamma)(y^+),\quad y^+\to\Phi(T\wedge
T_\kappa(\omega),\theta_{T\wedge T_\kappa(\omega)},\gamma)(y^+)
\]
is continuous. Hence by Castaing and Valadier \cite{CasVal77},
Lemma III.14, the above terms are measurable with respect to
$(\omega,y^+)$. The measurability follows now by the composition
formula (\ref{001}).
\end{proof}

\begin{remark}\label{r10}
{\rm i) Note that the solution of (\ref{eq11}), (\ref{eq11b}) can
be extended to any
time interval $[0,T]$. Then lemma \ref{l55}, \ref{l55} remain true for {\em any} $T>0$.\\
ii) Similar to the extension procedure we can show that
$\Psi(T,\omega,\gamma)$ is defined for any $T>0,\,\omega\in\Omega$
and $\gamma\in\mathcal{H}$. }
\end{remark}

%%%%%%%%%%%%%%%%%
%%%%%%%%%%%%%%%%%
%%%%%%%%%%%%%%%%%
%%%%%%%%%%%%%%%%%
\section{Existence of generalized fixed points} \label{ss}

By Theorem \ref{main},     the problem of finding
invariant manifolds  for a cocycle is  equivalent to finding
 generalized fixed points for a related  (but different) cocycle.
In this section, we present a generalized fixed point theorem
 for cocycles.

Let $\Omega$  and $\theta$ be  as in Section \ref{s2}, except that,
in this section, we do not need any measurability assumptions.
Namely, $\Omega$ is an invariant set (of full measure)
 under the metric dynamical system
 $\theta$.
Let $\Phi$ be a cocycle on a complete  metric space
$(\mathcal{G},d_{\mathcal{G}})$.

Recall that a mapping  $\gamma^\ast:\Omega\to   \mathcal{G} $
is called a {\em generalized
fixed point} of the cocycle $\Phi$ if
\[
\Phi(t,\omega,\gamma^\ast(\omega))=\gamma^\ast(\theta_t\omega)\;\text{for }t\in\mathbb{R}.
\]
Note that by the invariance of  $\Omega$ with respect to
$\{\theta_t\}_{t\in\mathbb{R}}$,  the trajectory $\mathbb{R}\ni
t\to\gamma^\ast(\theta_t\omega)\in \mathcal{G}$ forms an {\em entire
trajectory} for $\Phi$.\\

    The following  generalized fixed point theorem
for cocycles    is
similar to the third author's earlier work \cite{Schm97a}.

\begin{theorem}\label{t10}
Let $(\mathcal{G},d_{\mathcal{G}})$ be a  complete metric space with
 bounded metric.
Suppose that
\[
\Phi(t,\omega,\mathcal{G})\subset  \mathcal{G}
\]
for $\omega\in \Omega,\,t\ge 0$ and that $x\to\Phi(t,\omega,x)$ is
continuous. In addition, we assume the contraction condition: There
exists a constant $k<0$ such that for $\omega\in\Omega$
\begin{equation*}
%\label{equ21}
\sup_{x\not=y\in \mathcal{G}} \log
\frac{d_{\mathcal{G}} (\Phi(1,\omega,x),\Phi(1,\omega,y))} {d_{\mathcal{G}}(x,y)}\le k.
\end{equation*}
Then  $\Phi$ has a unique generalized fixed point $\gamma^\ast$ in
$\mathcal{G}$. Moreover, the following convergence property holds
\begin{equation*}
%\label{eq21a}
  \lim_{t\to\infty}\Phi(t,\theta_{-t}\omega,x)=\gamma^\ast(\omega)
\end{equation*}
for any $\omega\in\Omega$ and $x\in \mathcal{G}$.
\end{theorem}
\begin{proof}
Let $x\in  \mathcal{G}$. For $\omega\in \Omega$ we consider the sequence
\begin{equation}\label{u1}
n\to (\Phi(n,\theta_{-n}\omega,x)).
\end{equation}
 To see that this sequence is a Cauchy sequence, we compute  by using the cocycle property
\begin{align*}
d_\mathcal{G}(\Phi&(n,\theta_{-n}\omega,x),\Phi(n+1,\theta_{-n-1}\omega,x))\\
&= d_\mathcal{G}(\Phi(n,\theta_{-n}\omega,x),
\Phi(n,\theta_{-n}\omega,\Phi(1,\theta_{-n-1}\omega,x)))\\
&= d_\mathcal{G}(\Phi(1,\theta_{-1}\omega,\Phi(n-1,\theta_{-n}\omega,x)),
\Phi(1,\theta_{-1}\omega,\Phi(n-1,\theta_{-n}\omega,\Phi(1,\theta_{-n-1}\omega,x))))\\
&\le e^kd_\mathcal{G}(\Phi(n-1,\theta_{-n}\omega,x),
\Phi(n-1,\theta_{-n}\omega,\Phi(1,\theta_{-n-1}\omega,x)))\\
& \le e^{kn}d_\mathcal{G}(x,\Phi(1,\theta_{-n-1}\omega,x))
\end{align*}
for  $n\in\mathbb{N}$. We denote the limit of this Cauchy sequence by $\gamma^\ast(\omega)$.\\

If we replace  $x$ in (\ref{u1}) by another element $y\in  \mathcal{G}$ we
obtain the same limit which  follows from
\[
d_\mathcal{G}(\Phi(n,\theta_{-n}\omega,x),\Phi(n,\theta_{-n}\omega,y))\le
e^{kn}d_\mathcal{G}(x,y).
\]
This implies that $\gamma^\ast(\omega)$ is independent of choice
of $x$.
\\

Now we prove the convergence property
\[
\lim_{t\to\infty}\Phi(t,\theta_{-t}\omega,x)=\gamma^\ast(\omega).
\]
In fact,
\begin{align*}
d_\mathcal{G}&(\Phi(t,\theta_{-t}\omega,x),\Phi([t],\theta_{-[t]}\omega,x))\\
&=
d_\mathcal{G}(\Phi([t],\theta_{-[t]}\omega,\phi(t-[t],\theta_{-t}\omega,x)),\Phi([t],\theta_{-[t]}\omega,x))\\
&\le e^{k[t]}d_\mathcal{G}(\Phi(t-[t],\theta_{-t}\omega,x),x)\to
0\quad\text{for }t\to\infty
\end{align*}
where $[t]$ denotes the integer part of $t$. Since
$\Phi(t-[t],\theta_{-t}\omega,x)\in \mathcal{G}$ the values
$d_\mathcal{G}(\Phi(t-[t],\theta_{-t}\omega,x),x)$ are uniformly bounded for
$t\in\mathbb{R}$ and $x\in \mathcal{G}$.

Next, we show that $\gamma^\ast$ is, as a matter of fact,
a generalized fixed point for $\Phi$. Since $x\to\Phi(t,\omega,x)$ is
continuous, for $t\ge 0$ we obtain
\begin{align*}
\Phi&(t,\omega,\gamma^\ast(\omega))=\Phi(t,\omega,\lim_{n\to\infty}\Phi(n,\theta_{-n}\omega,x))\\
&=\lim_{n\to\infty}\Phi(t+n,\theta_{-n}\omega,x)=\lim_{n\to\infty}\Phi(t+n,\theta_{-n-t}\theta_t\omega,x)
=\gamma^\ast(\theta_t\omega).
\end{align*}

Finally, we prove the uniqueness of the generalized fixed point.
Suppose there is another generalized fixed point
$\bar\gamma^\ast(\omega)\in \mathcal{G}$. Let
$\Gamma^\ast=\{\gamma^\ast(\theta_t\omega),\,t\in\mathbb{R},\,\omega\in\Omega\}$
and $\bar \Gamma^\ast=\{\bar
\gamma^\ast(\theta_t\omega),\,t\in\mathbb{R},\,\omega\in\Omega\}$.
Since $\Gamma^\ast$ and $\bar \Gamma^\ast$ are bounded in $\mathcal{G}$ and
\begin{align*}
d_\mathcal{G}(\gamma^\ast(\omega),\bar\gamma^\ast(\omega))&=
d_\mathcal{G}(\Phi(n,\theta_{-n}\omega,\gamma^\ast(\theta_{-n} \omega)
),\Phi(n,\theta_{-n}\omega,\bar\gamma^\ast(\theta_{-n}
\omega)))\\
&\le e^{kn}\sup\{d_\mathcal{G}(x, y)\vert x\in \Gamma^\ast, y\in
\bar\Gamma^\ast\},
\end{align*}
letting $n\rightarrow \infty$, we have
$\gamma^\ast(\omega)=\bar\gamma^\ast(\omega)$. This completes the
proof.
\end{proof}

\begin{remark}
The constant  $ k$  in the above generalized fixed point theorem
may be taken as $\omega$-dependent,
as long as the following condition is satisfied:
\[
    \lim_{n\to\pm\infty}\frac{1}{n}\sum_{i=0}^{n-1}k(\theta_i\omega) =c<0.
\]
 This latter condition is usually assumed  in the situation of ergodicity.
 For applications see for instance Schmalfu{\ss} \cite{Schm97a}
 and Duan et al.
 \cite{DuaLuSchm02}.
\end{remark}

%%%%
%%%%%
%%%%%
%%%%%%%%%%%%%%
%%%%%%%%%%%%%%
\section{Random invariant manifolds}\label{s4}

In this final section, we show that  the random graph transform,
defined in (\ref{eq7}), has a generalized fixed point in the state space
\begin{equation}\label{space}
\mathcal{H}=B_{C_0^{0,1}}(0,\kappa)
 \;\; \mbox{with the metric}\;\;  d_\mathcal{H}(x,y):=\|x-y\|_{C_0^G},
\end{equation}
by using Theorem
\ref{t10}. Thus by
Theorem \ref{main}, the graph of this generalized fixed point
is an invariant manifold of the random dynamical system
generated by (\ref{u5}).

We first consider the basic properties of the
 metric space $\mathcal{H}$.

\begin{lemma}
The  metric space $\mathcal{H}=(B_{C_0^{0,1}}(0,\kappa),d_\mathcal{H}),\,d_\mathcal{H}(x,y):=\|x-y\|_{C_0^G}$ is
complete and   the metric $d_\mathcal{H}$ is bounded.
\end{lemma}
\begin{proof}
Let $(x_n)$ be a Cauchy sequence in $\mathcal{H}$. Since
$C_0^G(H^+;H^-)$ is complete we have $x_n\to x_0\in
C_0^G(H^+;H^-)$. Hence,  we have for any $y^+\in H^+$ that
$x_n(y^+)\to x_0(y^+)$. Subsequently,
\begin{equation}\label{u10}
\frac{\|x_n(y_1^+)-x_n(y_2^+)\|_H}{\|y_1^+-y_2^+\|_H}\to
\frac{\|x_0(y_1^+)-x_0(y_2^+)\|_H}{\|y_1^+-y_2^+\|_H}\quad\text{for
}n\to\infty
\end{equation}
for $y_1^+\not=y_2^+\in H^+$. Since the left hand side is
uniformly bounded by $\kappa$ so is the right hand side of
(\ref{u10}). Hence $x_0\in B_{C_0^{0,1}}(0,\kappa)$. The
boundedness assertion is easily seen.
\end{proof}

We now check the assumptions of the generalized  fixed point Theorem
\ref{t10}. Let $\Phi$ be the random dynamical system given by the
graph transform in (\ref{eq7}).

\begin{theorem}\label{t100}
Suppose that the gap condition (\ref{eq30}) is satisfied. Then the
 random  graph transform defined in (\ref{eq7})
has a unique generalized fixed point $\gamma^\ast(\omega,\cdot)$
in $\mathcal{H}$ where $\kappa$ is given in Lemma \ref{l20}. The
graph of this generalized fixed point, namely,
$M(\omega)=\{(x^+,\gamma^\ast(\omega,x^+)),\,x^+\in H^+\}$ is an
invariant manifold for the random dynamical  system $\phi$ generated
by (\ref{u5}).
\end{theorem}

\begin{proof}
By Lemma \ref{l20}, Theorem \ref{t11} we  know that $\Phi(T,\omega,\cdot)$ maps
$\mathcal{H}$
into itself.\\

Before we check the contraction condition in Theorem \ref{t10}
 we calculate an estimate
for $\|\Psi(1,\theta_1\omega,\gamma)\|_{C_0^G}$ for $\gamma\in
\mathcal{H}$. This norm is given by $\|w(0)\|_{C_0^G}$ where
$(w,v)$ is a solution of (\ref{eq8}) for $T=1$ and $\gamma\in
\mathcal{H}$. An estimate for $\|w(0)\|_{C_0^G}$ is given by
$W(0)$ defined in (\ref{eq11}), (\ref{eq11b}) for
$T=1,\,C=0,\,Y=1$. By the monotonicity of $W(0)$ in $\Gamma$ we
obtain that $W(0)$ for $\Gamma=\kappa=e_+^{-1}$ is an estimate of
$\|w(0)\|_{C_0^G}$ for {\em any} $\gamma\in\mathcal{H}$. Now we
can calculate $W(0)$ explicitly which gives us the estimate
\begin{equation}\label{u40}
\|w(0)\|_{C_0^G}\le
W(0)=e^{-\lambda_+-\int_0^1z(\theta_\tau\omega)d\tau}\quad (T=1!).
\end{equation}
We now check the contraction condition.
 To  this end we consider problem
(\ref{eq8}) for two different elements $\gamma_1,\,\gamma_2\in
\mathcal{H}$ and we denote the solutions  by
$w_i,\,v_i,\,i=1,\,2$. In particular, we have
\[
w_1(T)-w_2(T)=0,\quad
v_1(0)-v_2(0)=\gamma_1(w_1(0))-\gamma_2(w_2(0)).
\]
By  the Lipschitz continuity of  the nonlinear term $G$
in the random partial differential equation (\ref{u5}),  we can estimate
\[
\frac{\|\pi^\pm G(w_1+v_1)-\pi^\pm G(w_2+v_2)\|_H}{\|y^+\|_H} \le
L\frac{\|w_1-w_2\|_H+\|v_1-v_2\|_H}{\|y^+\|_H}
\]
which implies that
\[
\|\pi^\pm G(w_1+v_1)-\pi^\pm G(w_2+v_2)\|_{C_0^{G}} \le
L(\|w_1-w_2\|_{C_0^{G}}+\|v_1-v_2\|_{C_0^{G}}).
\]
Similar to Lemma \ref{l55} we  can estimate
\begin{equation}\label{eq202}
\|\Phi(1,\omega,\gamma_1)-\Phi(1,\omega,\gamma_2)\|_{C_0^G}=\|v_1(1)-v_2(1)\|_{C_0^G}
\end{equation}
by  $V(1)$ and $\|w_1(0)-w_2(0)\|_{C_0^G}$ by  $W(0)$, where
$V(t)$ and $W(t)$ is  a solution of (\ref{eq11}) with
\begin{equation}\label{equ17}
W(1)=0,\qquad V(0)=
\|\gamma_1-\gamma_2\|_{C_0^{G}}e^{-\lambda_+-\int_0^1z(\theta_\tau\omega)d\tau}+\kappa
W(0).
\end{equation}
Indeed, we can estimate the norm of initial condition
$v_1(0)-v_2(0)$:
\begin{align*}
\|v_1(0)-v_2(0)\|_{C_0^{G}}=&\\
\|\gamma_1(w_1(0))-\gamma_2(w_2(0))\|_{C_0^{G}} \le&
\|\gamma_1(w_1(0))-\gamma_2(w_1(0))\|_{C_0^{G}}\\
&+\|\gamma_2(w_1(0))-\gamma_2(w_2(0))\|_{C_0^{G}}\\
\le&\|\gamma_1-\gamma_2\|_{C_0^{G}}\|w_1(0)\|_{C_0^{G}}+\|\gamma_2\|_{C_0^{0,1}}\|w_1(0)-w_2(0)\|_{C_0^G}\\
\le&\|\gamma_1-\gamma_2\|_{C_0^{G}}e^{-\lambda_+-\int_0^1z(\theta_\tau\omega)d\tau}+
\kappa W(0).
\end{align*}
We have a bound for $\|w_1(0)\|_{C_0^{G}}$ from (\ref{u40}) and $\|\gamma\|_{C_0^{0,1}}\le\kappa$.
Since $V(1)$ as a solution (\ref{eq11}), (\ref{eq11b}) at $T=1$ is
increasing in $\Gamma$ and $C$ the value $V(1)$ for the above
generalized initial conditions (\ref{equ17}) is an estimate for
(\ref{eq202}) for any $\gamma_1,\,\gamma_2\in\mathcal{H}$. We have
chosen
$C=\|\gamma_1-\gamma_2\|_{C_0^{G}}e^{-\lambda_+-\int_0^1z(\theta_\tau\omega)d\tau}$.
\\

We now can calculate $V(1)$ explicitly. For these calculations we
have used that the solution operator $Q(t)$ for the linear problem
(\ref{eq11})  can be written as
\[
Q(t)[c_1,c_2] =c_1 \left(
\begin{array}{c}
e_+\\
1
\end{array}
\right) e^{\lambda_+ t+\int_0^tz(\theta_\tau\omega)d\tau}+ c_2
\left(
\begin{array}{c}
e_-\\
1
\end{array}
\right) e^{\lambda_- t+\int_0^tz(\theta_\tau\omega)d\tau}.
\]
These calculations of (\ref{eq11}) yield with the initial
conditions (\ref{equ17})
\begin{align*}
c_1&=e_-\|\gamma_1-\gamma_2\|_{C_0^{G}}\frac{-e^{-\lambda_+-\int_0^1z(\theta_\tau\omega)d\tau}}
{e_{+}-e_{-}}e^{\lambda_--\lambda_+}\\
c_2&=e_+\|\gamma_1-\gamma_2\|_{C_0^{G}}\frac{
e^{-\lambda_+-\int_0^1z(\theta_\tau\omega)d\tau}}{e_{+}-e_{-}}.
\end{align*}
In summary, we have for  $\gamma_1,\,\gamma_2\in\mathcal{H}$
\[
\|\Phi(1,\omega,\gamma_1)-\Phi(1,\omega,\gamma_2)\|_{C_0^G}=
\|v_1(1)-v_2(1)\|_{C_0^G}\le
V(1)=\|\gamma_1-\gamma_2\|_{C_0^{G}}e^{\lambda_--\lambda_+}.
\]
Since $\lambda_+>\lambda_-$,  we  thus obtain the contraction condition
in Theorem \ref{t10}
for $k=\lambda_--\lambda_+<0$.\\
We obtain similar estimates if we replace $T=1$ by $T>0$. Then
these estimates show us that
\[
\gamma\to\Phi(T,\omega,\gamma)
\]
is continuous at $\gamma\in \mathcal{H}$.\\

So  we have found that all assumption of Theorem \ref{t10}
are satisfied. Hence the dynamical system generated by the graph
transform  $\Phi$ has a unique generalized fixed point $\gamma^\ast$ in
$\mathcal{H}$. The graph of   $\gamma^\ast$   defines a desired
 invariant manifold for the random dynamical system
$\phi$  by Theorem \ref{main}.
\end{proof}

It remains to prove that this manifold is measurable.
\begin{lemma}
The manifold $M(\omega)$ is a random manifold.
\end{lemma}
\begin{proof}
The fixed point $\gamma^\ast(\omega,x^+)$ is the $\omega$-wise
limit of $\Phi(t,\theta_{-t}\omega,\gamma)(x^+)$ for $x^+\in H^+$
and for some $\gamma$ in $\mathcal{H}$ as  $t\to\infty$, see
Theorem \ref{t10}. Hence the mapping
$\omega\to\gamma^\ast(\omega,x^+)$ is measurable for any $x^+\in
H^+$. In order to see that $M$ is a random set we have to verify
that for any $x\in H$
\begin{equation}\label{eq2021}
\omega\to\inf_{y\in H}\|x-\pi^+y-\gamma^\ast(\omega,\pi^+y)\|_H
\end{equation}
is measurable, see Castaing and Valadier \cite{CasVal77} Theorem
III.9. Let $H_c$ be a countable dense set of the separable space
$H$. Then the right hand side of (\ref{eq2021}) is equal to
\begin{equation}\label{eq203}
\inf_{y\in H_c}\|x-\pi^+y-\gamma^\ast(\omega,\pi^+y)\|_H
\end{equation}
which follows immediately by the continuity of
$\gamma^\ast(\omega,\cdot)$. The measurability of (\ref{eq203})
follows since $\omega\to\gamma^\ast(\omega,\pi^+y)$ is measurable
for any $y\in H$.
\end{proof}
Under the additional assumption  $\hat\lambda>0>\check\lambda $ we
can show that $M$ is  an {\em unstable}  manifold denoted by
$M^+$: For any $\omega\in \Omega,\,t\ge 0$ and $x\in M^+(\omega)$
there exists an $x_{-t}\in M(\theta_{-t}\omega)$ such that
\begin{equation}\label{eq1000}
\phi(t,\theta_{-t}\omega,x_{-t})=x=x^++\gamma^\ast(\omega,x^+)
\end{equation}
and $x_{-t}$ tends to zero. We set
\[
x_{-t}=\Psi(t,\omega,\gamma^{\ast})(x^+)+
\gamma^\ast(\theta_{-t}\omega,\Psi(t,\omega,\gamma^\ast)(x^+)),\quad
x^+:=\pi^+x.
\]
Equation (\ref{eq1000})  is satisfied because
$x^+\to\pi^+\phi(t,\theta_{-t}\omega,x^++\gamma^\ast(x^+))$ is the
inverse of $x^+\to\Psi(t,\omega,\gamma^{\ast})(x^+)$, and because
$\gamma^\ast$ is the fixed point of the graph transform. The value
$\|\Psi(t,\theta_t\omega,\gamma^\ast)(x^+)\|_{C_0^G}$ can be
estimated by $W(0)$ a solution of (\ref{eq11}), (\ref{eq11b}) on
$[0,T]$ with $\Gamma=\kappa,\,C=0$ and $Y=1$ and
$\omega=\theta_{-t}\omega$. $W(0)$ can be  calculated explicitly
for any $T>0$. Hence
\begin{equation*}
%\label{equ31}
\|\Psi(t,\omega,\gamma^\ast)(x^+)\|_{H} \le
e^{-\lambda_+t-\int_{-t}^0z(\theta_\tau\omega)d\tau}\|x^+\|_H.
\end{equation*}
(We have to replace $\omega$ by $\theta_{-t}\omega$!) We can
derive from Lemma \ref{l100} iv)
\[
\int_{-t}^0 z(\theta_\tau\omega)d\tau<\eps t
\]
for any $\eps>0$ if $t$ is chosen sufficiently large depending on
$\omega$ and $\eps$. Hence
$\|\Psi(t,\omega,\gamma^\ast)(x^+)\|_{C_0^G}$ tends to zero
exponentially. On the other hand we have for $\gamma^\ast\in
\mathcal{H}$
\[
\|\gamma^\ast(\theta_{-t}\omega,\Psi(t,\omega,\gamma^\ast)(x^+))\|_{H}\le
\kappa \|\Psi(t,\omega,\gamma^\ast)(x^+)\|_{H}\to 0\quad \text{
for }t\to\infty.
\]
This convergence is exponentially fast. We conclude  that $M^+$ is
the unstable manifold for (\ref{u5}).
\\

However, our intention is to prove that (\ref{eq3a}) has an
invariant (unstable) manifold. On account of conjugacy of
(\ref{eq3a}) and (\ref{u5}) by (\ref{eq500}) and (\ref{eq501}) we
will now formulate the following result.

\begin{theorem} \label{spde}
Let $\phi$ by the random dynamical system generated by (\ref{u5})
and $\hat\phi$ be the solution version of (\ref{eq3a}) generated
by (\ref{eq503}). Then $M(\omega)$ is the invariant manifold of
$\phi$ if and only if $\hat M^+(\omega)=T^{-1}(\omega,M^+(\omega))$
is the invariant manifold of $\hat\phi$.  Moreover, if $M^+$ is an
unstable  manifold,  then so is $\hat M^+$.
\end{theorem}
\begin{proof}
We have  the relationship between $\phi$ and $\hat\phi$ given in
Lemma \ref{l111}
\begin{align*}
\hat\phi(t,\omega,&\hat M^+(\omega))=T^{-1}(\theta_t\omega,\phi(t,\omega,T(\omega,\hat M^+(\omega))))\\
&=T^{-1}(\theta_t\omega,\phi(t,\omega,M^+(\omega)))\subset
T^{-1}(\theta_t\omega,M^+(\theta_t\omega)) =\hat
M^+(\theta_t\omega).
\end{align*}
Note that $t\to z(\theta_t\omega)$ has a sublinear growth rate,
see Lemma \ref{l100}iii). Thus the transform
$T^{-1}(\theta_{-t}\omega)$ does not change the exponential
convergence of\\ $\Psi(t,\omega,\gamma^\ast(\omega))(x^+)$:
\[
\hat\Psi(t,\omega,\hat\gamma^\ast(\omega))
=T^{-1}(\theta_{-t}\omega,\Psi(t,\omega,T(\omega,\hat\gamma^\ast(\omega)))),\quad
\hat\gamma^\ast(\omega):=T^{-1}(\omega,\gamma^\ast(\omega)).
\]
It follows that $\hat M^+(\omega)$ is unstable.
\end{proof}

\begin{remark}
{\rm Note that  the main Theorem \ref{t100} represents the best possible
result in the following sense. If we consider the solution of the
two dimensional problem (\ref{eq11}) then this differential
equation generates a {\em non trivial} invariant manifold {\em if
and only if} the gap condition (\ref{eq30}) is satisfied. Hence we
can not formulate stronger {\em general} conditions for the
existence of global  manifolds. Here nontrivial means that the
dimension of the manifold is less than the dimension of the
space.}
\end{remark}
\appendix
\section{Proofs of the Lemmas \ref{l17}, \ref{l18} and \ref{l55}}
We now give the proof of the technical lemmas \ref{l17}, \ref{l18}
and \ref{l55}
which are based on the usual  Banach fixed point theorem.\\

{\bf Proof of Lemma \ref{l17}}:\\
We consider the following operator
\begin{align*}
\mathcal{T}_{T}&: C([0,T];C_0^{G}(H^+;H^+)\times C_0^{G}(H^+;H^-))\\
&\to C([0,T];C_0^{G}(H^+;H^+)\times C_0^{G}(H^+;H^-))
\end{align*}
for some $T>0$. Set $\mathcal{T}_{T}(w_1,v_1)=(w_2,v_2)$ where
\begin{align*}
\label{eq12}
\begin{split}
w_2(t)=&e^{\int_T^tz(\theta_r\omega)dr}\pi^+S(t-T)y^+\\
&-\int_t^{T}e^{\int_\tau^tz(\theta_r\omega)dr}\pi^+S(t-\tau)\pi^
+G(\theta_\tau\omega,w_1(\tau)+v_1(\tau))d\tau,\\
v_2(t)=&e^{\int_0^tz(\theta_r\omega)dr}\pi^-S(t)\gamma(w_2(0))\\
&+\int_0^{t}e^{\int_\tau^tz(\theta_r\omega)dr}\pi^-S(t-\tau)
\pi^-G(\theta_\tau\omega,w_1(\tau)+v_1(\tau))d\tau.
\end{split}
\end{align*}
Note that $w_1,\,v_1$ depend on $y^+,\,t,\,\omega$ and $\gamma$. A
fixed point for $\mathcal{T}_{T}$ is a solution of (\ref{eq8}) on
$[0,T]$. It is obvious that if
\[
(w_1,v_1)\in C([0,T];C_0^{G}(H^+;H^+)\times C_0^{G}(H^+;H^-))
\]
so is $(w_2,v_2)$. We check that the contraction condition of the
Banach fixed point theorem is satisfied. We set
\[
\Delta w_i=w_i-\bar w_i,\qquad \Delta v_i=v_i-\bar v_i,\quad
i=1,\,2.
\]
By the Lipschitz continuity of $\gamma$:
\[
\|\gamma(w_i(0))-\gamma(\bar w_i(0))\|_H\le L_\gamma\|\Delta
w_i(0)\|_H,\quad L_\gamma=\|\gamma\|_{C_0^{0,1}}.
\]
Hence we obtain by (\ref{eqvv}) for $H^+\ni y^+\not=0$
\begin{align*}
\frac{\|\Delta w_2(t)\|_H}{\|y^+\|_H} \le&
\int_t^{T}e^{\int_\tau^tz(\theta_r\omega)dr}e^{\hat\lambda(t-\tau)}L\frac{\|\Delta
w_1(\tau)\|_H+ \|\Delta v_1(\tau)\|_H}{\|y^+\|_H}d\tau
\\
\le&
L\int_t^{T}e^{\int_\tau^tz(\theta_r\omega)dr}e^{\hat\lambda(t-\tau)}d\tau
\bigg(\sup_{t\in [0,T]}\frac{\|\Delta w_1(t)\|_H}{\|y^+\|_H}
+\sup_{t\in [0,T]}\frac{\|\Delta v_1(t)\|_H}{\|y^+\|_H}\bigg)\\
\frac{\|\Delta v_2(t)\|_H}{\|y^+\|_H} \le&
L_\gamma\frac{\|\Delta w_2(0)\|_H}{\|y^+\|_H}e^{\int_0^tz(\theta_r\omega)dr}e^{\check\lambda t}\\
&+ \int_0^{t}e^{\int_\tau^tz(\theta_r\omega)dr}e^{\check\lambda
(t-\tau)}L\frac{\|\Delta w_1(\tau)\|_H+
\|\Delta v_1(\tau)\|_H}{\|y^+\|_H}d\tau\\
\le& L_\gamma
\int_0^{T}e^{\int_\tau^tz(\theta_r\omega)dr}e^{\hat\lambda(t-\tau)}L\frac{\|\Delta
w_1(\tau)\|_H+
\|\Delta v_1(\tau)\|_H}{\|y^+\|_H}d\tau\\
&+\int_0^{t}e^{\int_\tau^tz(\theta_r\omega)dr}e^{\check\lambda
(t-\tau)}L\frac{\|\Delta w_1(\tau)\|_H+
\|\Delta v_1(\tau)\|_H}{\|y^+\|_H}d\tau\\
\le& K(\omega,T,L_\gamma) \bigg(\sup_{t\in [0,T]}\frac{\|\Delta
w_1(t)\|_H}{\|y^+\|_H} +\sup_{t\in [0,T]}\frac{\|\Delta
v_1(t)\|_H}{\|y^+\|_H}\bigg).
\end{align*}
Choosing $T$ sufficiently small,  we have
\begin{equation}\label{000}
  K(\omega,T,L_\gamma)<1,\quad K(\omega,T,L_\gamma)=
LT\bigg((L_\gamma+1) e^{\int_0^T|z(\theta_r\omega)|+|\hat\lambda|
dr}d\tau+e^{\int_0^T|z(\theta_r\omega)| +|\check\lambda| dr}\bigg)
\end{equation}

We now can take the supremum with respect to $y^+\not=0$ and $t\in [0,T]$ for the left
hand side. Hence for
sufficiently small $T\le 1$ the operator $\mathcal{T}_{T}$ is a
contraction.
\hfill$\square$\\

{\bf Proof of Lemma \ref{l18}}:
\\
The proof of existence and uniqueness is similar to the proof in
Lemma \ref{l17}. The solution can be constructed by  successive
iterations of (\ref{eq11}), (\ref{eq11b}). If we start with
$V_1(t)\equiv\Gamma Y+C\ge\hat V_1(t)\equiv\hat\Gamma Y+\hat C$,
$\hat W_1(t)=W_1(t)\equiv Y$ we get
\[
V_2(t)\ge \hat V_2(t),\,W_2(t)\ge \hat W_2(t)\;,\cdots,\;
V_i(t)\ge \hat V_i(t),\,W_i(t)\ge \hat W_i(t)\;,\cdots.
\]
which gives the conclusion. These inequalities also show if
$(W(t),V(t))$ exist on $[0,T]$ so do $(\hat W(t),\hat V(t))$. The
inequalities for the contraction condition do not contain
$C$.\hfill$\square$
\\

{\bf Proof of Lemma \ref{l55}}:\\
Let $(w_i,v_i)$, $(W_i,V_i)$ be sequences generated by the
successive iterations starting with
$v_1(t)\equiv\gamma(y^+),\,w_1(t)\equiv y^+$ and $W_1=1,\,
V_1=L_\gamma=\|\gamma\|_{C_0^{0,1}}$. These sequences converge to the solution of
(\ref{eq8}) and (\ref{eq11})\,(\ref{eq11b}) provided $T$
sufficiently small. We then have
\begin{align*}
\|w_i(t)\|_{C_0^G}\le & e^{\int_T^tz(\theta_r\omega)+\hat\lambda
dr}
+\int^T_te^{\int_s^tz(\theta_r\omega)+\hat\lambda dr}\|\pi^+G(\theta_s\omega,w_{i-1}(s)+v_{i-1}(s))\|_{C_0^G}ds\\
\le &e^{\int^t_Tz(\theta_r\omega)+\hat\lambda dr}
+\int_t^Te^{\int_s^tz(\theta_r\omega)+\hat\lambda dr}(L\|w_{i-1}(s)\|_{C_0^G}+L\|v_{i-1}(s)\|_{C_0^G})ds\\
\|v_i(t)\|_{C_0^G}
\le & e^{\int_0^tz(\theta_r\omega)+\check\lambda dr}L_\gamma\|w_i(0)\|_{C_0^G}\\
&+\int_0^te^{\int_s^tz(\theta_r\omega)+\check\lambda
dr}(L\|w_{i-1}(s)\|_{C_0^G}+L\|v_{i-1}(s))\|_{C_0^G})ds
\end{align*}
and
\begin{align*}
W_i(t)=& e^{\int_T^tz(\theta_r\omega)+\hat\lambda dr}
+\int_t^Te^{\int_s^tz(\theta_r\omega)+\hat\lambda dr}
(L(W_{i-1}(s)+LV_{i-1}(s))ds\\
V_i(t)=& \L_\gamma
W_i(0)e^{\int_0^tz(\theta_r\omega)+\check\lambda dr}
+\int_0^te^{\int_s^tz(\theta_r\omega)+\check\lambda
dr}(LW_{i-1}(s)+LV_{i-1})ds.
\end{align*}
It is easily seen that $W_1(t)= \|w_1(t)\|_{C_0^G},\,V_1(t)\ge
\|v_1(t)\|_{C_0^G}$ and that if
\[
W_{i-1}(t)\ge \|w_{i-1}(t)\|_{C_0^G},\,V_{i-1}(t)\ge
\|v_{i-1}(t)\|_{C_0^G}
\]
then
\[
W_{i}(t)\ge \|w_{i}(t)\|_{C_0^G},\,V_{i}(t)\ge
\|v_{i}(t)\|_{C_0^G}
\]
which gives the conclusion. \hfill$\square$

%%%%%%%%%%%%%%%%%
%%%%%%%%%%%%%%%%%
%\section{Summary}

%We have considered invariant unstable manifolds for infinite
%dimensional non-autonomous dynamical systems, such as
%non-autonomous, or random, or stochastic, parabolic partial
%differential equations.

%By a random graphy transform and a random Banach
%contraction mapping fixed point theorem, we obtain
%  sufficient conditions, as in Theorems \ref{t100}  and \ref{spde},
%under which invariant unstable manifolds exist.

%%%%%%%%%%%%%%%%%
%%%%%%%%%%%%%%%%%


\begin{thebibliography}{1}

\bibitem{Arn98}
L.~Arnold.
\newblock {\em {R}andom {D}ynamical {S}ystems}.
\newblock Springer, New York, 1998.

\bibitem{BabVis89}
A.~B. Babin and M.~I. Vishik.
\newblock {\em Attractors of Evolution Equations}.
\newblock North-Holland, Amsterdam, London, New York, Tokyo, 1992.

\bibitem{BatLuZen98}
P.~Bates, K.~Lu, and C.~Zeng.
\newblock {\em Existence and Persistence of Invariant Manifolds for Semiflows
  in {B}anach Space}, volume 135 of {\em Memoirs of the AMS}.
\newblock 1998.

\bibitem{BenFla95}
A.~Bensoussan and F.~Flandoli.
\newblock Stochastic inertial manifold.
\newblock {\em Stochastics Stochastics Rep.}, 53(1--2):13--39, 1995.

\bibitem{CasVal77}
C.~Castaing and M.~Valadier.
\newblock {\em Convex Analysis and Measurable Multifunctions}.
\newblock LNM 580. Springer--Verlag, Berlin--Heidelberg--New York, 1977.

\bibitem{ChiLat97}
C.~Chicone and Y.~Latushkin.
\newblock Center manifolds for infinite dimensional non-autonomous differential
  equations.
\newblock {\em J. Diff. Eqns.}, 141:356--399, 1997.

\bibitem{ChoLuLin91}
S-N. Chow, K.~Lu, and X-B. Lin.
\newblock Smooth foliations for flows in banach space.
\newblock {\em Journal of Differential Equations}, 94:266--291, 1991.



\bibitem{DaPZab92}
G.~Da Prato and J.~Zabczyk.
\newblock {\em Stochastic Equations in Infinite Dimension}.
\newblock University Press, Cambridge, 1992.

\bibitem{DaPDeb96}
G.~Da Prato and A.~Debussche.
\newblock Construction of stochastic inertial manifolds using backward
  integration.
\newblock {\em Stochastics Stochastics Rep.}, 59(3--4):305--324, 1996.

\bibitem{DuaLuSchm02}
J.~Duan, K.~Lu, and B.~Schmalfu{\ss}.
\newblock Unstable manifolds for equations with time dependent coefficients.
\newblock 2002.
\newblock In preparation.

\bibitem{GirChu95}
T.~V. Girya and I.~D. Chueshov.
\newblock Inertial manifolds and stationary measures for stochastically
  perturbed dissipative dynamical systems.
\newblock {\em Sb. Math.}, 186(1):29--45, 1995.

\bibitem{Had01}
J.~Hadamard.
\newblock Sur l'iteration et les solutions asymptotiques des equations
  differentielles.
\newblock {\em Bull. Soc. Math. France}, 29:224--228, 1901.

\bibitem{Hen81}
D.~Henry.
\newblock {\em Geometric theory of semilinear parabolic equations}, volume 840
  of {\em Lecture Notes in Mathematics}.
\newblock Springer-Verlag, New York, 1981.

\bibitem{KokSie01}
N.~Koksch and S.~Siegmund.
\newblock Pullback attracting inertial manifolds for nonautonomous dynamical
  systems.
\newblock {\em J. Dyn. Differ. Equations}, 2002.
\newblock To appear.

\bibitem{Lia47}
A.~M. Liapunov.
\newblock {\em Probl\`eme g\'eneral de la stabilit\'e du mouvement}, volume~17
  of {\em Annals Math. Studies}.
\newblock Princeton, N.J, 1947.

\bibitem{Kun90}
H.~Kunita.
\newblock {\em Stochastic Flows and Stochastic Differential Equations}.
\newblock Cambridge University Press, Cambridge, 1990.

\bibitem{MohScheu99}
S.-E.~A. Mohammed and M.~K.~R. Scheutzow.
\newblock The stable manifold theorem for stochastic differential equations.
\newblock {\em The Annals of Probability}, 27(2):615--652, 1999.

\bibitem{Oks92}
B.~{\O}ksendale.
\newblock {\em Stochastic Differential Equations}.
\newblock Springer--Verlag, Berlin--Heidelberg--New York, third edition, 1992.

\bibitem{Per28}
O.~Perron.
\newblock {\"U}ber {S}tabilit\"at und asymptotisches {V}erhalten der
  {I}ntegrale von {D}ifferentialglei\-chungssystemen.
\newblock {\em Math. Z.}, 29:129--160, 1928.

\bibitem{Rue82}
D.~Ruelle.
\newblock Characteristic exponents and invariant manifolds in Hilbert spaces.
\newblock {\em Ann. of Math.}, 115:243--290, 1982.

\bibitem{Schm97c}
B.~Schmalfu{\ss}.
\newblock The random attractor of the stochastic {L}orenz system.
\newblock {\em ZAMP}, 48:951--975, 1997.

\bibitem{Schm97a}
B.~Schmalfu{\ss}.
\newblock A random fixed point theorem and the random graph transformation.
\newblock {\em Journal of Mathematical Analysis and Applications},
  225(1):91--113, 1998.

\bibitem{Schm99a}
B.~Schmalfu{\ss}.
\newblock Attractors for the non-autonomous dynamical systems.
\newblock In K.~Gr{\"o}ger,  B.~Fiedler and J.~Sprekels, editors, {\em
  Proceedings {EQUADIFF99}}, pages 684--690. World Scientific, 2000.

\bibitem{Sel67}
G.~R. Sell.
\newblock Non-autonomous differential equations and dynamical systems.
\newblock {\em Amer. Math. Soc.}, 127:241--283, 1967.

\bibitem{CarLanRob01}
T.~Caraballo, J.~Langa  and J.~C. Robinson.
\newblock A stochastic pitchfork bifurcation in a reaction-diffusion equation.
\newblock  2001.
\newblock Submitted.

\bibitem{Vis92}
M.~I. Vishik.
\newblock {\em Asymptotic Behaviour of Solutions of Evolutionary Equations}.
\newblock {C}ambridge University Press, {C}ambridge, 1992.

\bibitem{Wanner} T. Wanner.
\newblock Linearization random dynamical systems.
\newblock  In C. Jones, U. Kirchgraber and H. O. Walther, editors,
{\em Dynamics Reported}, Vol. 4, 203-269, Springer-Verlag, New
York, 1995.


\end{thebibliography}
\end{document}